\documentstyle[12pt]{article}
\newtheorem{thm}{Theorem}
\newtheorem{dfn}{Definiton}
\newtheorem{lem}{Lemma}
\newtheorem{prop}{Proposition}
\newtheorem{cor}{Corollary}
\newcommand{\To}{{\bf T}}
\newcommand{\Ho}{H\"{o}lder}
\newcommand{\bmu}{\mbox{\boldmath $\mu$}}
\newcommand{\T}{Teichm\"{u}ller}
\title{The \T \ space of an Anosov diffeomorphism of $T^{2}$}
\author{Elise E. Cawley}
\date{July 17, 1991}

\begin{document}

\maketitle

\section{Introduction}

In this paper we consider the space of smooth conjugacy classes of an Anosov
diffeomorphism of the two-torus.  
A diffeomorphism $f$ of a manifold $M$ is Anosov if there is a
continuous invariant splitting of the tangent bundle 
$TM = E^{s} \bigoplus E^{u}$, where
the subbundle $E^{s}$ is contracted by $f$, and $E^{u}$ is expanded.  More
precisely, if $\| \cdot \|$ is a Riemannian  metric on $M$, there
are constants $c > 0$ and $\lambda < 1$ such that
\begin{eqnarray*}
\| Df^{n} \cdot v \| & \leq &  c \lambda^{n} \| v \|  
                       \ \mbox{for v in} 
                       \  E^{s} \\
\| Df^{-n} \cdot v \| & \leq &  c \lambda^{n} \| v \|  
                       \ \mbox{for v in} 
                       \  E^{u} 
\end{eqnarray*}

\noindent for all positive integers $n$.  If $f$ is $C^{1 + \alpha}$ for
$0 < \alpha < 1$, it can be shown that the splitting is in fact H\"{o}lder
continuous.

The only
2-manifold that supports an Anosov diffeomorphism is the 2-torus \cite{F}.
Franks and Manning showed that  every Anosov diffeomorphism of $\To^{2}$ is
topologically conjugate to a linear example; that is, to an 
automorphism defined by a hyperbolic element of $GL(2,{\bf Z})$ whose 
determinant has absolute value one \cite{F}, \cite{Ma}.  Consider
$f$ and $g$ which are topologically conjugate, so there is a homeomorphism
taking the orbits of $f$ to the orbits of $g$.  If the conjugacy is in fact
smooth, then $f$ and $g$ must have the same expanding and contracting 
eigenvalues at corresponding periodic points.  De la Llave, Marco, and 
Moriyon have shown that the eigenvalues at periodic points are a {\em
complete} smooth invariant:  if the eigenvalues of $f$ are the same as
the eigenvalues of $g$ at periodic points that correspond under a 
topological conjugacy, then the conjugacy is smooth \cite{L},\cite{MM1},
\cite{MM2}.

The question  arises:
what sets  of eigenvalues occur as the Anosov diffeomorphism ranges over a 
topological conjugacy class?  The information in the set of expanding 
eigenvalues of $f$ is recorded by a H\"{o}lder cyclic cohomology class
associated  
to $f$.  The real-valued function
\[x \mapsto \phi_{u}(x)  = -{\rm log} 
\parallel Df(x) \parallel _{u}\;,\] where 
$\parallel Df \parallel _{u}$ is the Jacobian of $f$ along the unstable 
bundle $E^{u}$, 
can be described as a ``cocycle'' over $f$.  
The cohomology class of this cocycle,
that is, its residue class modulo the space of coboundaries $x \mapsto
u(f(x)) - u(x)$,
is independent of the choice of Riemannian metric on $M$.  Moreover, the
cohomology class of a H\"{o}lder cocycle defined over an Anosov diffeomorphism
is determined by the sums of values of the cocycle over the various 
periodic orbits, which  
for $\phi _{u}$ is simply minus the logarithm of the expanding eigenvalue at
the periodic point. Similarly, the information in the set of contracting
eigenvalues is recorded by the cohomology class of the cocycle defined 
by $\phi_{s} = {\rm log} \parallel
Df \parallel _{s}$ where $\parallel Df \parallel _{s}$ is the Jacobian of
$f$ along the stable bundle $E^{s}$.  The sign convention that makes these
cocycles negative is chosen for consistency  with the notation of Bowen,
Ruelle, and Sinai in the theory of Gibbs and equilibrium measures for Anosov
systems.  

The question asked above can be reformulated: what pairs of cohomology 
classes (one determined by the expanding eigenvalues, and one by the 
contracting
eigenvalues) occur as the diffeomorphism ranges over a topological conjugacy
class?  The cohomology is defined over the entire conjugacy class by pulling
the Jacobian cocycles back to a fixed representative of the conjugacy class.  
The purpose
of this paper is to answer this question:  {\em all} pairs of H\"{o}lder
reduced  cohomology classes occur. (The reduced cohomology is the cyclic
cohomology 
divided out by the constant cocycles.  The pair of reduced cohomology classes
is still sufficient information to determine the smooth conjugacy class.)

The Teichm\"{u}ller space $T(f)$ of an Anosov diffeomorphism $f$ 
is defined to be the set of smooth structures 
preserved by the topological dynamics
determined by $f$.  This is the smooth category version of the 
Teichm\"{u}ller space of a
rational map, which was studied by McMullen and Sullivan \cite{MS}.
We show that for an Anosov diffeomorphism of $\To^{2}$, there is a 
natural bijection from the  Teichm\"{u}ller
space to the  product $G(f) \times G(f^{-1})$, where $G(f)$ denotes
the real-valued \Ho \ reduced cyclic cohomology
over $f$.

\noindent The main theorem is:

\begin{thm}
Let $f:{\bf T}^{2} \rightarrow {\bf T}^{2}$ be an Anosov diffeomorphism.  
Then there is a natural bijection  
$T^{1 +  H}(f) \leftrightarrow 
G(f) \times G(f^{-1})$.
\end{thm}

\noindent An easy corollary will be

\begin{thm}
Let $f:{\bf T}^{2} \rightarrow {\bf T}^{2}$ be a volume preserving Anosov
diffeomorphism. 
Then there is a natural bijection $T_{vol}^{1 +  H}(f) \leftrightarrow  G(f)$.
\end{thm}

\noindent  Here
$G(f)$ is the H\"{o}lder reduced cyclic cohomology over $f$, where the
H\"{o}lder  
exponent is allowed to vary between $0$ and $1$. More precisely:  let 
$C^{\alpha}$ denote the space of $\alpha$-\Ho \ functions on $T^2$.
Let $C^{H} = \cup_{\alpha \in (0,1)} C^{\alpha}$.  Then $G(f)$ is 
the quotient of $C^{H}$ by the subspace of ``almost coboundaries'' 
\[x \mapsto u(f(x))  - u(x) + K\;,\] where $u$ is a function on $T^2$ and 
$K \in {\bf R}$.  (It follows in this setting  that $u \in C^{H}$).
$T^{1 +  H}(f)$ is
the Teichm\"{u}ller space of $C^{1+H}$ invariant smooth structures.
An Anosov diffeomorphism is called {\em volume preserving} if it admits
an invariant measure that is absolutely continuous with repect to 
Lebesgue.  $T_{vol}^{1 + H}
(f)$ is the restriction of this Teichm\"{u}ller space to the volume preserving
elements. See the appropriate sections for precise definitions.

The main theorem can be restated as follows:

\bigskip
\noindent {\bf \large Theorem $1'$}
{\em Let $L$ be a hyperbolic automorphism of the torus $\To^{2}$.  Given two
\Ho \
functions $\phi_{u}$ and $\phi_{s}$ from $\To^{2}$ to $\bf R$,  there exist 
uniquely defined constants $P_{u}$ and $P_{s}$, and a unique $C^{1 + H}$
smooth
structure on $\To^2$ which is preserved by $L$,  and determines the
cohomology classes
$< \phi_{u} - P_{u}>$ and $<\phi_{s} - P_{s}>$.
Moreover, $L$ is Anosov in this smooth structure.}

\bigskip

\noindent {\bf Remark.}  The \Ho \ exponent of the new
smooth structure depends on the \Ho \ norms and exponents of the pair of
functions, 
and on a dynamically defined  norm (the {\em Bowen}, 
or {\em variation} norm) of the cohomology classes. \cite{B3}

\bigskip

\noindent {\bf Remark.}  The topological conjugacy between 
$C^{1 + H}$ Anosov diffeomorphisms is in fact \Ho \ continuous \cite{Mn}.
Therefore  this description
is independent of the choice of ``base point'' as the linear mapping.

\begin{cor}
Let $\lambda_{u}$ and $\lambda_{s}$ denote the unstable and stable 
eigenvalues at a periodic point of period $n$.  The numbers
$|\lambda_{u}|^{1/n}$ and $|\lambda_{s}|^{1/n}$ 
can be prescribed arbitrarily on any finite set 
of periodic points, up to (non-unique) constant factors ${\rm exp}(-P_{u})$
and ${\rm exp}(P_{s})$, respectively.
\end{cor}

\noindent {\bf Remark.}  The correction factors in the corollary are 
asymptotically unique along any sequence of periodic point sets $P_i$ with
the 
property that the normalized dirac mass on $P_{i}$ converges to 
Haar measure (which is also the measure of maximal entropy) on ${\bf
T}^{2}$.
\bigskip

We give a sketch of the proof of Theorem 1.  The main step is to show
that a  cohomology class over an
Anosov diffeomorphism $f$ determines  a canonical invariant transverse
measure class to 
the stable foliation. 
This is constructed using the Gibbs measure class defined by the cohomology
class.
When the stable foliation is co-dimension 1, the
transverse measure class can be interpreted as a transverse smooth
structure.  In dimension 2,
when both foliations are co-dimension 1, the two transverse structures
(determined by
the two given cohomology classes) define 
a product smooth structure, which is invariant by $f$.

Ruelle and Sullivan \cite{RS}, and Sinai \cite{Si}, gave a  transverse
interpretation of 
a particular Gibbs measure, namely that associated to the constant cocycle.
In 
this case, and only this case, one has a transverse measure, as opposed to a
transverse
measure class.  Under the isomorphism of Theorem 1, this corresponds to the
linear
member of the topological conjugacy class.  The present paper shows how to
extend the 
decomposition of a Gibbs measure into transverse stable and transverse
unstable part,
as described in \cite{RS} for the constant cocycle measure, to the general
case.  

The organization of the paper is as follows.  In section 2 we
define the Teichm\"{u}ller space $T(f)$.  In section 3 we recall facts
about cyclic cohomology over a ${\bf Z}$ action,  over  a foliation or an
equivalence relation, and  describe the Jacobian and
Radon-Nykodym cocycles. In Section 4, the map which gives the isomorphism in
Theorem 1 is described.  Section 5 describes the cocycle properties of 
Gibbs measures, and collects some necessary results. 
Section 6 gives the statements of  the transverse structure 
realization results, and proves Theorem 1 assuming these.  The main body of
the
proof is in Section 7, where the tranverse measure class is constructed.
Section 8 gives a simple description of the smooth structure defined by a 
pair of cohomology classes in terms of explicit coordinates on rectangles in
a 
Markov partition. 

\noindent {\bf Acknowledgements.}  It is a pleasure to thank Dennis Sullivan
for
many hours of conversation, for inspiration and encouragement.  
Jack Milnor made many helpful comments on the manuscript.  
I also thank IHES and the Institute for Mathematical Sciences at 
Stony Brook for 
their hospitality while this work was being completed.

\section{The Teichm\"{u}ller space of an Anosov diffeomorphism}

Let $f:M \rightarrow M$ be a $C^{r}$ Anosov diffeomorphism, $0 < r \leq
\omega$.
The $C^{r}$ \T \ space $T^{r}(f)$ is defined as follows.  Consider triples
$(h, N, g)$ where $g:N \rightarrow N$ is a $C^{r}$ Anosov diffeomorphism,
and $h:M \rightarrow N$ is a homeomorphism satisfying $g \circ h = h \circ
f$.
We call such a triple a {\em marked Anosov diffeomorphism modeled on 
$f:M \rightarrow M$}.  Two such triples $(h_{1},N_{1},g_{1})$ and
$(h_{2},N_{2},g_{2})$ are {\em equivalent} if the
homeomorphism  $s:N_{1} \rightarrow N_{2}$ defined by 
$s \circ h_{1} = h_{2}$ is in fact a $C^{1}$ diffeomorphism.
The \T \ space $T^{r}(f)$ is the space of equivalence classes of 
triples. 

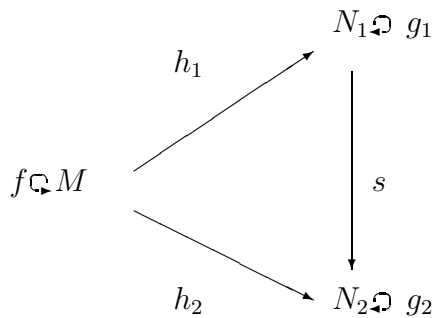
\begin{figure}[t]
\caption{\T \ space of $f:M \rightarrow M$.}
\begin{center}
\begin{picture}(195,165)

\put(13,90){$f$}
\put(30,90){$M$}
\put(150,90){$s$}
\put(75,135){$h_{1}$}
\put(75,45){$h_{2}$}
\put(135,150){$N_{1}$}
\put(135,45){$N_{2}$}
\put(163,150){$g_{1}$}
\put(163,45){$g_{2}$}

\put(60,97.5){\vector(3,2){68}}

\put(60,82.5){\vector(2,-1){68}}
\put(142.5,135){\vector(0,-1){75}}

\put(24,93){\oval(6,6)[l]}
\put(24,93){\oval(6,6)[t]}
\put(29,90){\vector(1,0){0}}

\put(154,153){\oval(6,6)[r]}
\put(154,153){\oval(6,6)[t]}
\put(149,150){\vector(-1,0){0}}

\put(154,48){\oval(6,6)[r]}
\put(154,48){\oval(6,6)[t]} 
\put(149,45){\vector(-1,0){0}}

\end{picture}
\end{center}
\end{figure}

We also define the $C^{1 + H}$ \T \ space $T^{1 + H}(f)$.  Let 
 $f:M \rightarrow M$ be a $C^{1+\alpha}$ diffeomorphism for some $0 < 
\alpha < 1$.  Consider all marked Anosov diffeomorphisms $g:N \rightarrow
N$ modeled on $f$ where $g$ is $C^{1 + \alpha^{\prime}}$ for some\break 
$0 < \alpha^{\prime}  < 1$ (that is, the \Ho \ exponent of $g$
is not necessarily the same as that of $f$).  Then $T^{1 + H}$ is the space
of 
equivalence classes, where equivalence is defined just as it was in the 
$C^{r}$ category.

\section{Cyclic cohomology.}

\subsection{Group actions and the Jacobian and Radon-Nykodym cocycles.}

We describe here the notion of cocycle over a group action, and 
the associated notions of coboundary and cohomology (see \cite{Z} and 
\cite{K}.)  We work in the topological category since all cocycles we
are interested in have at least this degree of regularity.
We consider $\Gamma$, a locally compact, second countable group,
and a continuous right action of $\Gamma$ on a topological space $M$.  
We will be especially interested
in the case of a ${\bf Z}$ action defined by a diffeomorphism of a 
manifold $M$.  

A real-valued (additive) cocycle over the action of $\Gamma$ is a
continuous map
\[
\Phi:M \times \Gamma \rightarrow {\bf R}
\]
\noindent satisfying:
\[
\Phi(x,\gamma_{1}\cdot\gamma_{2}) = \Phi(x,\gamma_{1}) + 
                                     \Phi(x \cdot \gamma_{1}, \gamma_{2})
\]

\noindent Here $x \rightarrow x \cdot \gamma$ denotes the
action of $\gamma$ on the point $x$. A {\em coboundary} is a cocycle of 
the form $\Phi(x,\gamma) =
u(x\cdot \gamma) - u(x)$ where $u:M \rightarrow {\bf R}$ is a
continous function. 
The function $u$ is called the
{\em transfer function} of the coboundary $\Phi$.  
Two cocycles are {\em equivalent} or {\em cohomologous} if their 
difference is a coboundary.
The cohomology over
the action of $\Gamma$ is the space of equivalence classes of cocycles.
A cocycle of the form $\Phi(x,\gamma) = u(x\cdot \gamma) - u(x) +
K(\gamma)$, where
$K:\Gamma \rightarrow {\bf R}$ is a homomorphism,
 is called an {\em almost coboundary}.  The space of cocycles
modulo almost coboundaries is the {\em reduced cohomology} over the action.
If $\Phi$ is a cocycle, we denote its cohomology class and reduced
cohomology
class $<\Phi>$ and $<\Phi>_{*}$, respectively.
When the space $M$ has additional structure, we can consider e.g. \Ho,
Lipshitz, or smooth cocycles (and coboundaries and cohomology).

The cohomology equivalence relation on cocycles has the following meaning.
Let $F(M,{\bf R})$ be the space of continuous maps from $M$ to ${\bf R}$.
Then a cocycle $\Phi$ defines an action of $\Gamma$ on $F(M,{\bf R})$
as follows.  If $T:M \rightarrow {\bf R}$, then $(\gamma \cdot T)(x) =
\Phi(x,\gamma) + T(x \cdot \gamma)$.  Note that if $\Phi$ is the identically
$0$ cocycle, then this action is just the usual pull-back of functions by
the group action.  When the 
two cocycles $\Phi$ and $\Psi$  are cohomologous, the actions they define
are equivalent.   If $u$ is the transfer function of the coboundary 
that relates
$\Phi$ to $\Psi$, then the map $U:F(M,{\bf R}) \rightarrow F(M,{\bf R})$
defined by $T \rightarrow T + u$ is an isomorphism which conjugates the 
action defined by $\Phi$ to the action defined by $\Psi$.
\bigskip

\noindent {\bf Example 1: The Jacobian cocycle.}  Suppose that $\Gamma$ 
acts by diffeomorphisms on a Riemannian manifold $M$.  We define the
(additive) Jacobian cocycle $J:M  \times \Gamma \rightarrow {\bf R}$
by $J(x,\gamma) = {\rm log} \parallel D\gamma (x) \parallel$ where
$\parallel
\cdot \parallel$ is the Riemannian metric.  The chain rule for 
differentiation is 
precisely the cocycle condition.  If we choose a new Riemannian metric
$\parallel \cdot \parallel _{1}$
on $M$, then the new Jacobian cocycle is cohomologous to the original one.
The transfer function is simply the logarithm of the ratio of volume 
elements with repect 
to the two metrics.  
Hence there is a cohomology class, the Jacobian class $<J>$,  naturally 
associated to a smooth
group action on a smooth manifold.  The apriori coarser invariant, the
reduced Jacobian class $<J>_{*}$, is also defined.
\bigskip

\noindent {\bf Example 2: The Radon-Nykodym cocycle.}
Suppose $\Gamma$ acts on the measure space $(M,\mu)$, quasi-preserving
the measure $\mu$.  The (additive) Radon-Nykodym cocycle 
$R:M \times \Gamma \rightarrow
{\bf R}$ is defined by $R(x,\gamma) = {\rm log}
\frac{d\mu(\gamma(x))}{d\mu(x)}$.
(Since we have retricted to the topological category, we assume that 
the Radon-Nykodym is continuous.)  Again the cocycle condition is
the chain rule.  If $\nu$ is a measure that is equivalent to $\mu$, with 
Radon-Nykodym derivative $r = \frac{d\nu}{d\mu}$, then the corresponding
cocycles are cohomologous, via the transfer function $u = {\rm log (r)}$.
So  there is a cohomology class, the Radon-Nykodym class $<R>$, and a 
reduced cohomology class $<R>_{*}$,  naturally 
associated to an action
that preserves a measure class. 
\bigskip

We now focus on the case $\Gamma = {\bf Z}$.  A cocycle over a ${\bf Z}$
action is determined by its values on the generator: $\phi(x) =: \Phi(x,1)$.
The cocycle condition implies that $\Phi(x,n) = \sum_{k=0}^{n-1} \phi(x\cdot
k)$ 
where $x \cdot k$ denotes the action of $k$ on the point $x$.  If 
$f:M \rightarrow M$ is the action of the generator, then we write
$\Phi(x,n) = \sum_{k=0}^{n-1} \phi \circ f^{k}(x)$.  Hence the cocycles over
a ${\bf Z}$ action can be identifed with continuous functions on $M$.  A
coboundary is a function of the form $u \circ f - u$.  An almost coboundary
is a function of the form $u \circ f - u + K$, where $K$ is a constant.

\bigskip

\noindent{\bf Example 3: the BRS classes of an Anosov diffeomorphism.}
Suppose $f:M \rightarrow M$ is an Anosov diffeomorphism, with $TM = E^{s}
\bigoplus E^{u}$.  Let $\parallel \cdot \parallel$ be a Riemannian metric
on $M$.  Since the subbundle $E^{u}$ is preserved by $f$, we can define the 
{\em unstable Jacobian cocycle} to be the cocycle over the action of $f$ 
determined by the function $\phi_{u}(x) = -{\rm log} \parallel Df(x)
\parallel _{u}$.
Here $\parallel Df(x) \parallel _{u}$ denotes the Jacobian in  the unstable
direction.  (The minus sign is a convention in the theory of Bowen, Ruelle, 
and Sinai.)
Similarly
we can define a stable Jacobian cocycle $\phi_{s}(x) = {\rm log} \parallel
Df(x)
\parallel_{s}$.  (With this sign convention, the stable Jacobian cocycle
of $f$ is the same as the unstable Jacobian cocycle of $f^{-1}$.)  If 
$g$ is an Anosov diffeomorphism which is smoothly conjugate to $f$, say
$g \circ h = h \circ f$, then the corresponding Jacobian cocycles are
cohomologous.  Namely, $\phi_{u}^{g} \circ h  = \phi_{u}^{f} + u \circ f -
u$ where
$u = -{\rm log} \parallel Dh \parallel_{u}$, and $\phi_{s}^{g} \circ h =
\phi_{s}^{f} +
v \circ f - v$ where $v = {\rm log} \parallel Dh \parallel_{s}$.  In other
words, the smooth conjugacy class of $f$ naturally determines a pair of 
cohomology classes, $(<\phi_{u}>,<\phi_{s}>)$, which we will refer to
as the unstable and stable {\em BRS classes} (for Bowen, Ruelle, and Sinai)
of $f$.
The corresponding reduced classes will be referred to as the reduced BRS
classes of $f$.

\subsection{Cocycles over a foliation.}

It will be useful in what follows to have the notion of a cocycle over a
foliation
$\cal F$ of a manifold $M$.  The definition depends on the {\em graph}, or 
holonomy groupoid, of the foliation, which was constructed by Winkelnkemper
\cite{Co},\cite{W}.  
An element $\gamma$ of the graph $GR({\cal F})$ is a pair of points of $M$,
$x = s(\gamma)$ and $y = r(\gamma)$, together with an equivalence class of
smooth
paths $\gamma(t)$ from $x$ to $y$, tangent to the foliation.  
Two paths $\gamma_{1}(t)$, $\gamma_{2}(t)$ from $x$ to $y$ are equivalent if
the 
holonomy of the path $\gamma_{2} \circ \gamma_{1}^{-1}$ from $x$ to $x$ is
the identity.
There is a natural composition law on $GR$, defined when the endpoint of one
pair
is the first point of another pair. A  cocycle over the foliation $\cal F$
is a function
\[
\Phi:GR({\cal F}) \rightarrow {\bf R}
\] 
such that
\[
\Phi(\gamma_{1} \circ \gamma_{2}) = \Phi(\gamma_{1}) + \Phi(\gamma_{2})
\]
whenever the composition $\gamma_{1} \circ \gamma_{2}$ is defined.
A  coboundary is a cocycle of the form $\Phi(\gamma) = u(r(\gamma)) -
u(s(\gamma))$
where $u$ is a function on $M$.  
If the leaves of $\cal F$ have trivial holonomy, then the point $\gamma \in
GR({\cal F})$
depends only on the pair $x = r(\gamma)$ and $y = s(\gamma)$.  In this case
we will
use the notation $\Phi(x \rightarrow y)$  for the cocycle. 
\bigskip

\subsection{Cocycles over an equivalence relation}

Let ${\cal F} \subset M \times M$ be defined by an equivalence relation 
$\sim$ on $M$.
Then a cocycle over ${\cal F}$ is a function
\[
\Phi : {\cal F} \rightarrow {\bf R}
\]
satisfying 
\[ 
\Phi(x,z) = \Phi(x,y) + \Phi(y,z)
\]
whenever $y$ and $z$ belong to the equivalence class of $x$. See \cite{HK}  
A coboundary is a cocycle of the form $\Phi(x,y) = u(y) - u(x)$ where 
$u: M \rightarrow {\bf R}$.

\section{The Bowen-Ruelle-Sinai map}

Let $f:\To^{2} \rightarrow \To^{2}$ be a $C^{1 + \alpha}$ Anosov
diffeomorphism.  Let $G(f)$
denote the \Ho \ reduced  
cohomology over $f$, where the \Ho \ exponent is allowed to vary. 
The Bowen-Ruelle-Sinai map,
\[
BRS:T^{1 + H}(f) \rightarrow G(f) \times G(f^{-1})
\] 
\noindent is defined as follows.  Let $(h,N,g)$ be a representative of a
point in the \T \ space
$T^{1 + H}(f)$.  Let $\phi_{u}(g)$ and $\phi_{s}(g) = \phi_{u}(g^{-1})$ be
the 
unstable and stable Jacobian cocycles of $g$.  We map the \T \ point to the
pair of reduced 
cohomology classes $(<\phi_{u}(g) \circ h >_{*}, < \phi_{u}(g^{-1}) \circ
h>_{*})$.
Because the topological conjugacy between two $C^{1 + \alpha}$ 
Anosov diffeomorphisms is
always \Ho, these are \Ho \ cohomology classes.  
The image point is independent of the choice of representative, because a 
smooth conjugacy changes the Jacobian cocycles by a coboundary.

\begin{thm}[de la Llave, Marco, Moriyon]
The map $BRS$ is injective.
\end{thm}

\noindent {\bf Remark.}  Two points need to be added to the theorem proved
by 
de la Llave, Marco, and Moriyon (\cite{L},\cite{MM1},\cite{MM2})  
to give the stated result.  They consider 
$C^{r}$ diffeomorphisms where $2 \leq r \leq \omega$.  They prove that 
the pair of BRS cohomology classes is a complete $C^{r}$ conjugacy
invariant.  The first 
point is that if $\phi$ and $\psi$ are BRS cocycles over an Anosov
diffeomorphism
whose difference is an almost coboundary, then in fact the difference is a
coboundary. (See the remarks about the {\em pressure} in the next section.)
The second point is that the $C^{1 + H}$ smoothness
case is simpler than the higher smoothness case, because the foliations have
$C^{1 + H}$ transverse smoothness.  An easy modification of the ideas of
de la Llave, Marco, and Moriyon
proves that the BRS cohomology classes are a complete $C^{1 + H}$
conjugacy invariant.
\bigskip

To prove Theorem 1, we need to show that the map $BRS$ is {\em surjective}.
That is, given a pair 
of reduced cohomology classes $(<\phi_{u}>_{*}, <\phi_{s}>_{*})$ over $f$, 
we need to construct a marked
Anosov diffeomorphism $(h,N,g)$ modeled on $f$, whose pair of 
unstable and stable reduced BRS cohomology classes,
pulled back to $f:M  \rightarrow M$ is this pair.

This problem can be restated as follows.  We consider the smooth torus 
$M$
to be defined by a $C^{1 + \alpha}$ system of charts $(U_{\beta},
\eta_{\beta})$  on a topological torus $\To^{2}$.  The mapping
$f$ is defined on the topological torus, and is assumed to define a $C^{1 +
\alpha}$
Anosov diffeomorphism when viewed in the smooth charts
$(U_{\beta},\eta_{\beta})$.  
Given a pair of reduced cohomology classes $(<
\phi_{u}>_{*},<\phi_{s}>_{*})$ over
$f$, the problem is to construct a new smooth system of charts
$(U_{\beta}^{\prime},
\eta_{\beta}^{\prime})$ on the topological torus $\To^{2}$, with the
property that
in these charts the mapping $f$ is $C^{1 + \alpha^{\prime}}$, for some
$0 < \alpha^{\prime} < 1$, and the reduced BRS cohomology classes in this
smooth
structure are $<\phi_{u}>_{*}$ and $<\phi_{s}>_{*}$.  In other words, we
identify
the underlying point sets of the smooth tori $M$ and $N$ via the
homeomorphism
$h$, and vary the point in \T \ space by varying the smooth structure on
this
point set.
\bigskip

\noindent {\bf Proof of Theorem 2 from Theorem 1.} An Anosov diffeomorphism 
is volume preserving if and only if the forward and backward BRS cohomology
classes coincide (under the canonical identification of the forward and 
backward cohomology) \cite{B}.  Hence the BRS map restricted to the volume
preserving
diffeomorphisms maps onto the ``diagonal'' in $G(f) \times G(f^{-1})$,
which is naturally identified with $G(f)$.

\section{Gibbs measures and associated cocycles}

In this section we recall some properties of Gibbs measures.  We emphasize
the dynamically defined {\em equivalence relations} and associated 
cocycles that define a 
Gibbs measure.  The transverse measure class constructed in this paper
can be viewed as a ``transverse Gibbs measure.''  It is obtained by
focusing on a different dynamically defined equivalence relation and 
its associated cocycle.

Gibbs measures are defined for dynamical systems with a local product 
structure.  See Baladi's thesis \cite{Ba} for a nice exposition and 
proofs of some basic results for the general case.  Here we only need
consider the two classes:  Anosov diffeomorphisms and subshifts of finite
type.  The latter arise because an Anosov diffeomorphism has a 
presentation as a quotient of a subshift of finite type. 

Let $A$ be an $r \times r$ matrix of $0$'s and $1$'s.  Consider the setof
bi-infinite sequences $\Sigma_{A} = \lbrace {\bf x} \rbrace$,
where ${\bf x} = \ldots x_{-k}x_{-k+1}\ldots x_{-1}x_{0}x_{1}
\ldots x_{l}x_{l+1}\ldots$
is in $\Sigma_{A}$ if and only if $x_{i} \in \lbrace 1,\ldots,r \rbrace$ and
$A_{x_{i}x_{i+l}} = 1$ for all $i$. The shift map $\sigma:\Sigma_{A}
\rightarrow
\Sigma_{A}$ is defined by $\sigma({\bf x}) = {\bf y}$ where $y_{i} =
x_{i+1}$.
We give $\Sigma_{A}$ the topology defined by the metric
$d({\bf x},{\bf y}) = \Sigma_{x_{i} \neq y_{i}}
2^{-| i |}$.
The metric space $\Sigma_{A}$ together with the shift map is called the
{\em subshift of finite type} defined by $A$.

If ${\bf x} \in \Sigma_{A}$, we define the {\em stable set} $W^{s}({\bf x})$
to
be the set of all ${\bf y} \in \Sigma_{A}$ such that $d(\sigma^{n}
{\bf x},\sigma^{n}{\bf y}) \rightarrow 0$ as $n \rightarrow \infty$.
Similarly we
define the {\em unstable set} $W^{u}({\bf x})$ to be the set of ${\bf y} \in
\Sigma_{A}$ such that $d(\sigma^{-n}{\bf x},\sigma^{-n}{\bf y}) \rightarrow
0$
as $n \rightarrow \infty$.
A {\em local stable set} $W^{s}_{\epsilon}({\bf x})$
is the set of points
${\bf y}$ such that $d(\sigma^{n}({\bf x}),\sigma^{n}({\bf y})) < \epsilon$
for $n \geq 0$.
Local unstable sets are defined similarly.  A stable set has an intrinsic
topology defined by the metric $d^{-}({\bf x},{\bf y}) =
\sum_{x_{i} \neq y_{i}} 2^{i}$.  Similarly, an
unstable set has an intrinsic topology defined by the metric 
$d^{+}({\bf x},{\bf y}) = \sum_{x_{i} \neq y_{i}} 2^{-i}$.  We will refer
to the collection of stable (unstable) sets as the stable (unstable) 
foliation of $\Sigma_{A}$.  

If $f: M \rightarrow M$ is an Anosov diffeomorphism, the stable foliation
${\cal W}^{s}$ and unstable foliation ${\cal W}^{u}$ are 
defined similarly,
and have tangent distributions $E^s$ and $E^u$, respectively.

Anosov diffeomorphisms and subshifts of finite type have a {\em
local product structure}. For every point $x$ in $M$ (or 
$\Sigma_{A}$) there is a neighborhood $U$ of $x$, and a 
homeomorphism
\[
u: W^{s}_{\epsilon}(x) \times W^{u}_{\epsilon}(x) \rightarrow U
\]
that takes verticals $\lbrace w \rbrace \times W^{u}_{\epsilon}(x)$
onto local unstable sets, and horizontals  $W^{s}_{\epsilon}(x) \times 
\lbrace z \rbrace$ onto local stable sets.

There is an equivalence relation on $M$, and $\Sigma_{A}$, defined
by the pair of foliations ${\cal W}^{s}$ and ${\cal W}^{u}$. Namely,
\[
x \sim y \Leftrightarrow x \in W^{s}(y) \cap W^{u}(y).
\]
If $x \sim y$, then there are neighborhoods $U_{x}$ of $x$ and $U_{y}$
of $y$, and a homeomorphism
\[
\theta : U_{x} \rightarrow U_{y}
\]
such that $z \sim \theta (z)$ for all $z \in U_{x}$.  These are called
{\em conjugating homeomorphims}.  The pseudogroup of conjugating
homeomorphisms generates this equivalence relation (referred to in the
sequel as the Gibbs equivalence relation).

\bigskip
\noindent The following definition is due to Capocaccia \cite{Ca}.
\begin{dfn}
Let $f: M \rightarrow M$ be Anosov.  Let $\phi : M \rightarrow {\bf R}$
be continuous.  A measure $\mu$ on $M$ is a {\em Gibbs measure for $\phi$}
if 
\[
{\rm log}(\frac{d\mu(\theta(x))}{d\mu(x)}) = \sum_{k = - \infty}^{\infty}
			      (\phi \circ f^{k} ( \theta(x)) 
			      	- \phi \circ f^{k}(x))
\]
for every conjugating homeomorphism $\theta$.
\end{dfn}
It is implicit in the definition that both sides of the equation are 
well-defined.
The left hand side of the equation is  the logarithmic
Radon-Nykodym cocycle associated to the Gibbs equivalence relation
and the measure $\mu$.  The function $\phi$ should be regarded as
a cocycle over $f$.  If $\phi$ is changed by an almost coboundary
$u \circ f - u + K$ the expression on the right hand side does not
change.  Hence a Gibbs measure is associated to a reduced cohomology
class over $f$, and is defined by an associated cocycle (referred 
to as the Gibbs cocycle) over the 
Gibbs equivalence relation.  The definition of Gibbs measure 
in the subshift of finite type case is exactly analagous.

Now we describe the cocycle properties of the transverse measure 
class we are going to construct. Instead of the Gibbs equivalence
relation, we consider only the stable foliation, and the corresponding
holonomy pseudogroup.  If $\phi : M \rightarrow {\bf R}$
is continuous, and if 
\begin{equation}
\Phi(x \rightarrow y) =: \sum_{k = 0}^{\infty}
			(\phi \circ f^{k} (y)
                                - \phi \circ f^{k}(x))
\end{equation}
is finite whenever $x \in W^{s}(y)$, and this expression 
defines a continuous function
of $x$ on a small transversal to ${\cal W}^{s}$, then $\Phi$ is a cocycle
over $W^{s}$, which we will refer to as the transverse Gibbs 
cocycle.  Moreover, if $\phi$ is changed by an almost coboundary
$u \circ f - u + K$, then $\Phi$ changes by the coboundary
$U(x \rightarrow y) = -u(y) + u(x)$.  The transverse Gibbs measure
class to be constructed will have the following properties.  We 
associate to a reduced cohomology class $< \phi >_{*}$ over the 
mapping the cohomology class $< \Phi >$ over the stable foliation
defined by equation 5.1.  If $\mu$ is a representative measure on 
a small transversal $\tau$, then the logarithmic Radon-Nykodym cocycle
of $\mu$ over the holonomy ${\rm hol}: \tau \rightarrow \tau$ will be of the
form
$\Phi(x \rightarrow y) + u(y) - u(x)$, where $u: \tau \rightarrow {\bf R}$
is a continuous function. (In fact $u$ will be \Ho \, but we postpone
the discussion of this until later.)

\bigskip

\noindent The following result of Bowen leads us to consider \Ho \ cocycles
over
$f$.

\begin{prop}[Bowen] \cite{B}  
If $\phi : M \rightarrow {\bf R}$ is \Ho, 
then the associated Gibbs cocycle and transverse Gibbs cocycle
exist, and are \Ho.  
\end{prop}
That is, the Gibbs and transverse Gibbs cocycles are \Ho \ on the 
domains of the homeomorphisms in the associated pseudogroups.
(In fact they are \Ho \ cocycles over the relevant {\em metric}
equivalence relations, but we will not need this.)
The proposition is true for \Ho \ cocycles over a subshift of 
finite type, as well.

\noindent {\bf Remark.}  There is a larger class of cocycles for 
which the Gibbs cocycles are finite.  In fact there is a natural 
norm identified by Bowen (the {\em variation norm}) which defines
a Banach space of cocycles with well-defined associated Gibbs
cocycles. \cite{B3}  One easily carries out the construction of a
transverse 
{\em continuous} measure class for these cocycles. 
An open problem is to determine
a regularity description  of these transverse measure classes.

The \Ho \
cohomology  over an Anosov diffemorphism is naturally associated to the
topological conjugacy class of the map.  This is
because the conjugacy between two $C^{1 + H}$ Anosov diffeomorphisms
is always \Ho \ continuous. 
The \Ho \ coboundaries over an Anosov diffemorphism form a {\em closed}
subspace
of the \Ho \ cocycles.  This is
a consequence of the Livshitz theorem, which states that the cohomology
class
of a \Ho \ cocycle $\phi$  over an Anosov diffeomorphism 
is determined by the values of
the cocycle over periodic orbits, i.e. by the sums $\sum_{k=0}^{n-1} \phi
\circ
f^{k}(p)$ where $f^{n}(p) = p$ \cite{Li}.  Since these sums 
are $0$ for a coboundary, and
this is a closed condition, the coboundaries are a closed subspace.

\bigskip

We now collect various results which will be needed in the sequel.  
The {\em one-sided subshift of finite type} $\Sigma_{A}^{+}$
associated to  a $0-1$ 
matrix $A$ is defined just as in the subshift of finite type case,
but one considers only one-sided sequences.  The shift map 
$\sigma : \Sigma_{A}^{+}  \rightarrow \Sigma_{A}^{+}$ is then
an expanding endomorphism.

\noindent The following lemma is due to Sinai and Bowen \cite{B}:

\begin{lem}
Let $\phi: \Sigma_{A} \rightarrow {\bf R}$ be a \Ho \ function.  Then $\phi$
is cohomologous (via a \Ho \ transfer function)
to a \Ho \ function $\phi^{+}$ with the property that $\phi^{+}(\bf x)$
depends only on the forward part of the sequence, 
i.e. on $x_{0},x_{1},x_{2},\ldots$.
\end{lem}

\noindent So the \Ho \ cyclic cohomologies over $(\Sigma_{A}, \sigma)$
and over $(\Sigma_{A}^{+}, \sigma)$ are isomorphic. 

We give a brief description of Bowen's proof, as we will need to
know the form of $\phi^{+}$.
Let $[i]$ denote the ``rectangle'' consisting of those  sequences
${\bf x}  \in \Sigma_{A}$  with  $x_{0} = i$.
Let $W^{u}({\bf x},i)$ be the local unstable set through ${\bf x}$
intersected with $[i]$.
For each $i$, choose ${\bf x}^{i} \in [i]$.   Define $r:\Sigma_{A}
\rightarrow
\Sigma_{A}$ by projecting along local stable sets
in the rectangle $[i]$, onto $W^{u}({\bf x}^{i},i)$.
Let $u:\Sigma_{A} \rightarrow {\bf R}$ be defined by
\[
u({\bf y}) = \sum_{k=0}^{\infty}(\phi \circ \sigma^{k}({\bf y}) 
- \phi \circ \sigma^{k}
(r({\bf y})))
\]
The function $u$ is \Ho \ along local unstable sets, by  Proposition 1
above in the subshift of finite type setting.  It can be checked that
\[
\phi^{+} = \phi + u \circ \sigma - u
\]
depends only on the forward part of a sequence.

\begin{flushright}
$\Box$
\end{flushright}

In the general setting of a homeomorphism $f:M \rightarrow M$ of a compact
metric space, the {\em pressure} function $P:C(M) \rightarrow {\bf R}$
is defined on the space $C(M)$ of continuous functions on $M$. 
See \cite{Wa3} or \cite{B}.   In fact, the
pressure is defined on the set of cohomology classes:  if $\phi$ is
cohomologous to $\psi$, then $P(\phi) = P(\psi)$.  In addition,
$P(\phi + K) = P(\phi) + K$.  
For the purposes of the present paper, the pressure can be viewed as 
defining an imbedding of the reduced cohomology into the
(unreduced) cohomology.  Namely, the reduced class $<\phi>_{*}$ is 
mapped to the unreduced class $<\phi - P(\phi)>$. The image is 
precisely the set of cohomology classes with pressure zero.
We note that if $f$ is an Anosov diffeomorphism, and
$<\phi>$ is the unstable BRS class, then $P(<\phi>) = 0$ \cite{B}.

\begin{thm}[Bowen-Ruelle-Sinai Theorem]
Let $\Sigma_{A}^{+}$ be a transitive one-sided subshift of finite type 
(transitive means that there is a dense orbit).  Let
$\phi : \Sigma_{A}^{+} \rightarrow {\bf R}$ be a \Ho \ function with
$P(\phi) = P$.
Then there is a  unique probability measure $\mu$ on
$\Sigma_{A}^{+}$ such that ${\rm log}\frac
{d\mu(\sigma({\bf x}))}{d \mu({\bf x})} = -\phi({\bf x}) + P$.
\end{thm}

\noindent See \cite{B},\cite{Wa2}. The measure $\mu$ is positive on open
sets.

The measure $\mu$ has the following ``one-sided Gibbs'' property.
Suppose that for some $n > 0$, we have $\sigma^{n}({\bf x})
= \sigma^{n}({\bf y})$.  Then there is a homeomorphism  $T_{{\bf x}{\bf y}}$
from a neighborhood of ${\bf x}$ to a neighborhood of ${\bf y}$,
defined by the property $\sigma^{n}(T_{{\bf x}{\bf y}}({\bf z}))
= \sigma^{n}({\bf z})$.
For $\mu$ as in the BRS Theorem, we have
\[
\log(\frac{d\mu(T_{{\bf x}{\bf y}}({\bf z}))}{d\mu({\bf z})}) =
\sum_{k = 0}^{n-1}(\phi (\sigma^{k}(T_{{\bf x}{\bf y}}({\bf z})))
- \phi (\sigma^{k} ({\bf z}))).
\]

\noindent We will need the following 

\begin{cor}[Local Uniqueness Corollary]
Let $\Sigma_{A}^{+}$ be a transitive one-sided subshift of finite type.  Let
$\phi : \Sigma_{A}^{+} \rightarrow {\bf R}$ be a \Ho \ function.  
Let $P = P(\phi)$.   Let $\nu$ be
a finite measure, supported  on an open set $V \subset \Sigma_{A}^{+}$.
Suppose
that $\nu$ satisfies the property: ${\rm log}(\frac{d\nu(\sigma^{n}({\bf
x})}
{d\nu({\bf x})}) = -\sum_{k=0}^{n-1}\phi \circ \sigma^{k}({\bf x}) + nP$,
whenever
${\bf x }\in V$ and $\sigma^{n}({\bf x}) \in V$.
Then $\nu$ coincides, up to a constant factor, with the
measure $\mu$ associated to $\phi$ by  
the Bowen-Ruelle-Sinai Theorem, restricted to $V$.
\end{cor}

\noindent {\bf Proof of Corollary 1.}  There is an open subset $U \subset V$
and an $n$ such that $\sigma^{n}$ is injective on $U$ 
and $\sigma^{n}(U) = \Sigma_{A}^{+}$.
Define a measure $\tilde{\nu}$ on all of $\Sigma_{A}^{+}$
by
\[
\frac{d\tilde{\nu}(\sigma^{n}({\bf x}))}{d\nu({\bf x})} =
{\rm exp}(-\sum_{k=0}^{n-1}\phi \circ \sigma^{k}({\bf x}) + nP)
\]
where ${\bf x} \in U$.
Then $\tilde{\nu}$ agrees with $\nu$ on $V$, by the derivative hypothesis on
$\nu$,
and ${\rm log}\frac{d\tilde{\nu}(\sigma({\bf x}))}{d\tilde {\nu}({\bf x}}
= - \phi({\bf x}) + P$.  Now
apply the uniqueness part of the Bowen-Ruelle-Sinai Theorem.

\begin{flushright}
$\Box$
\end{flushright}

\begin{prop}[The invariant measure.]
Let $\Sigma_{A}^{+}$ be a transitive subshift of finite type, and let
$\phi:\Sigma_{A}^{+} \rightarrow {\bf R}$ be \Ho.  Then there exists
a \Ho \ function $h$ with the following property. Let $\phi^{\prime} =
\phi + h - h \circ \sigma - P(\phi)$.  Then $\phi^{\prime} < 0$ and
$\mu_{\phi^{\prime}}$ (the measure associated to 
$ \phi^{\prime}$ by the BRS theorem) is invariant under $\sigma$.
\end{prop}

\noindent See \cite{Le},\cite{Wa}.
Let $\mu_{\phi}$ and $\mu_{\phi^{\prime}}$ be the measures associated
to $\phi$ and $\phi^{\prime}$, respectively, 
by the BRS theorem.  Then the construction yields
\[
\frac{d \mu_{\phi^{\prime}}}{d \mu_{\phi}} = \log ( h).
\]
Since $P(\phi^{\prime}) = 0$, the Radon-Nykodym derivative of
$\mu_{\phi^{\prime}}$
under the shift map is $\exp(-\phi^{\prime}) \geq  c > 1$.  Hence the
invariant  measure is
expanded by the shift map.  

\bigskip

We now address an important subtlety of the smooth structure
construction.  Ultimately we obtain a smooth structure
by integrating a representative measure on a transversal.  A different 
representative measure differs by a \Ho \ Radon-Nykodym
derivative, that is {\em \Ho \ with respect to the underlying
metric of, say, the linear toral diffeomorphism.}  We need 
to know that the Radon-Nykodym derivative is \Ho \ with 
respect to the {\em new} smooth coordinate, namely the measure
itself.  This will follow from the following proposition.

Suppose we have a mixing one-sided subshift of finite
type $\Sigma_{A}^{+}$.
A cylinder set $C_{n}$ of length $n$ associated to a finite word $x_{0}x_{1}
\ldots x_{n-1}$ is the set of all sequences in $\Sigma_{A}^{+}$
that begin with this word.  Suppose we have a \Ho \  map $\pi:\Sigma_{A}^{+}
\rightarrow I$ onto  an interval $I \subset {\bf R}$ satisfying the
following properties.
First,  $\pi$ is such that  $\pi({\bf x}) = \pi({\bf y})$
implies that $\pi(\sigma({\bf x})) = \pi(\sigma({\bf y}))$.
Second, we assume that the image by $\pi$ of a cylinder set is
an interval.  Let $\phi^{\prime} $ be a \Ho \ function on I.
Let $\mu$ be the measure on $\Sigma_{A}^{+}$, associated to the pull-back of
$\phi$
by $\pi$, constructed
in the Bowen-Ruelle-Sinai theorem.
Let $\mu_{0}$ be the measure corresponding to the constant function.
Assume that  $\pi$ is injective on a set of full
measure with respect to both $\mu_{0}$ and $\mu$.
There are two metrics $d_{0}$ and $d_{\phi}$ defined on $I$ by the
push-forward by $\pi$ of $\mu_{0}$ and $\mu$, respectively.

\begin{prop}
The identity map $\iota: (I,d_{0}) \rightarrow (I,d_{\phi})$ is
quasisymmetric.
\end{prop}

\noindent See \cite{Ja} and \cite{Ji} for the proof of quasisymmetry.
We give the outline of the proof in the appendix.

\begin{cor}
A function on $I$ is \Ho  \ in the $d_{0}$ metric if and only if
it is \Ho \ in the $d_{\phi}$ metric.
\end{cor}

\noindent {\bf Proof.}  A quasisymmetric homeomorphism is \Ho \ \cite{A}.

\section{Realizing cohomology classes as transverse structures}

\subsection{Transverse measure class to a foliation}

Let $M$ be a smooth $n$-dimensional manifold, and let ${\cal F}$ be a
$k$-dimensional foliation
of $M$. That is, $M = \cup_{F \in {\cal F}} F$ where each $F \in {\cal F}$
is a smooth
submanifold of $M$.  $M$ is covered by {\em flow-boxes} $D^{k} \times
D^{n-k}$ with the 
property that each leaf $F \in {\cal F}$ meets a flow-box in a collection of
disks of the
form $D^{k} \times \lbrace y \rbrace$.  A {\em transversal} $\tau$ to the
foliation is a 
smooth $(n-k)$-dimensional submanifold that meets each leaf $F$
transversely.  

A {\em transverse measure} $\mu_{\cal F}$ assigns to each small transversal
a measure, 
with finite total
mass, with the property that the measure is invariant under the holonomy
pseudogroup
\cite{Co},\cite{RS}.

By relaxing the condition on invariance under holonomy, we arrive at the
notion of a 
{\em transverse measure class}, $\bmu_{\cal F}$.  This object assigns to
each small
transversal  a  measure {\em class}, with finite total mass, which is
invariant under 
the holonomy pseudogroup.  If the foliation has sufficient  transverse
regularity so that 
it preserves a smoothness class $\Lambda$, e.g.
where $\Lambda$ denotes \Ho, Lipshitz, or $C^{r}$ regularity, then we can
define 
a {\em transverse $\Lambda$ measure class} by requiring that the
representative measures
on a transversal are equivalent with Radon-Nykodym derivatives in the class
$\Lambda$.

\bigskip

\noindent {\bf Example.}  Suppose $\cal F$ is a foliation of $M$ with
transverse 
regularity $C^{ 1 + H}$, i.e. the holonomy maps on transversals are 
$C^{1 + H}$.  Then Lebesgue measure on transversals defines a transverse
\Ho \  measure class.

\bigskip

Let $f:M \rightarrow M$ be a diffeomorphism, preserving the foliation $\cal
F$.
Then we say that $f$ preserves the transverse $\Lambda$ measure class $\bmu$
if for every
small transversal $\tau$, and measurable subset $E \subset \tau$, $E$ has
positive
$\bmu$-measure if and only if $f(E) \subset f(\tau)$ has positive
$\bmu$-measure,
and moreover the Radon-Nykodyn derivative has regularity $\Lambda$.

There is a $\Lambda$ cohomology class over $f$ naturally associated to an
$f$-invariant
transverse $\Lambda$ measure class $\bmu_{\cal F}$.  It can be defined as
follows.  
Pick a finite
covering of $M$ by flow-boxes $B_{i} = D^{k}_{i} \times D^{n-k}_{i}$.  
Choose
representative measures $\mu_{i}$ on transversals $\tau_{i} = \lbrace x_{i}
\rbrace
\times D^{n-k}_{i}$ in each flow-box.  
Let $\lbrace \alpha_{i} \rbrace$ be a smooth partition of unity subordinate
to the 
covering by flow-boxes.  Define 
\[
\phi(x) = -{\rm log} \frac{d\mu(f(x))}{d\mu(x)}
\]
where $\mu = \sum_{i} \alpha_{i}(x) \mu_{i}(x)$,
regarded as a measure on the local quotient space obtained by projecting
along
the leaf factors $D_{i}^{k} \times \lbrace y \rbrace$ in the flow boxes
containing
$x$.

The  $\Lambda$ cohomology class of $\phi$ is independent of the choice
of covering by flow-boxes, the  representative
measures, and the partition of unity.  
We will call this the Radon-Nykodym class of $f$ acting on the 
transverse measure class $\bmu_{\cal F}$.

In a similar way, we can define the notion of a {\em transverse smooth
structure}.
This is an assignment of a smooth structure to each small transversal, with
the 
property that the holonomy pseudogroup acts smoothly with respect to this
smooth
structure.  Note that the assigned smooth structure on a transversal in 
general will have 
nothing to do with that induced on a transversal by the ambient 
smooth structure of the manifold $M$.  If, as above, $f$ is a diffeomorphism
preserving the foliation $\cal F$, then an {\em f-invariant transverse 
smooth structure} is one in which the action of $f$ on transversals is
smooth
(with respect to the assigned structure).
There is a natural cohomology class associated to the action of $f$, defined
as in the transverse measure class case, but with Radon-Nykodym derivative
replaced by the Jacobian.  We will refer to this as the transverse
Jacobian class of the action of $f$ on the transverse smooth structure.

\subsection{Radon-Nykodym realization}

We are now ready to describe the main step in the construction of an
invariant smooth
structure from a pair of cohomology classes over an Anosov diffeomorphism
$f$.

\begin{thm}[Radon-Nykodym realization.]
Let $f:M \rightarrow M$ be an Anosov diffeomorphism.  Let $<\phi>_{*}$ be a
\Ho \
reduced cohomology class over $f$.  Then there is an f-invariant transverse
\Ho \ measure
class $\bmu$ to the stable foliation $W^{s}$, with the property  
that the reduced Radon-Nykodym class of $f$ acting on $\bmu$ is
$<\phi>_{*}$. 
\end{thm}

The transverse measure class $\bmu$ in the theorem has the additional
property that
the measure class on a transversal is positive on open subsets of the
transversal.  When
$M = T^{2}$, transversals to $W^{s}$ are one-dimensional.  In this case, an
$\alpha-$
\Ho \ transverse measure
class that is positive on open sets is equivalent to a $C^{1 +  \alpha}$
transverse 
smooth 
structure. 

\begin{thm}[Jacobian realization in dimension 2]
Let $f:T^{2} \rightarrow T^{2}$ be Anosov.  Let $<\phi>_{*}$ be a \Ho \
reduced
cohomology class
over $f$.  Then there is an f-invariant transverse $C^{1 +  \alpha}$ smooth
structure to the stable foliation $W^{s}$, with the property that the 
reduced Jacobian class of $f$ acting on this transverse smooth structure is 
$<\phi>_{*}$.

\end{thm}

\noindent The Radon-Nykodym realization theorem is proved in section 7.

\subsection{Complementary transverse smooth structures}

Let $\cal F$ and $\cal G$ be foliations of complementary dimension.
We say that the foliations  {\em intersect transversely} if
the leaves of $\cal F$ and $\cal G$ meet transversely.  We assume that there
is a system of simultaneous flow-boxes of the form 
\[
D^{k} \times D^{l}
\]
where each leaf of $\cal F$ meets a flow-box in a collection of disks of the
form $D^{k} \times \lbrace y \rbrace$, and each leaf of $\cal G$ meets a 
flow-box in a collection  of disks of the form $\lbrace x \rbrace \times
D^{l}$.

Let $\cal F$ and $\cal G$ be foliations of complementary dimension,
intersecting
transversely.  Then a pair of transverse smooth structures, one for each
foliation,
determines a canonical smooth structure on $M$ as follows.  
Consider a simultaneous flow-box $D^{k} \times D^{l}$.  Pick a point $(x,y)$
in the flow-box.  Then the disk $\lbrace x \rbrace \times D^{l}$ is a
transversal
to the foliation $\cal G$.  Similarly the disk $D^{k} \times \lbrace y
\rbrace$
is a transversal to the foliation $\cal F$. So there is a product smooth
structure
on the flow-box determined by the transverse smooth structures on these
disks.
The overlap maps for these charts are block diagonal, with the blocks being 
the derivative of the holonomy for each foliation.  So the product structure
has the same degree of smoothness as the transverse structures.
\bigskip

\noindent {\bf Proof of Theorem 1} 
We show that the BRS map is surjective,
assuming the Jacobian realization theorem. Let $<\phi_{u}>_{*}$ and 
$<\phi_{s}>_{*}$ be the reduced cohomology classes which we want to realize.
We have complementary, transverse foliations
$W^{s}$ and $W^{u}$, which by the Jacobian realization lemma can be equipped
with a transverse smooth structures, with associated reduced Jacobian
classes 
equal to $<\phi_{u}>_{*}$ and $<\phi_{s}>_{*}$.  We 
define a product structure as just described. In this smooth structure,
the unstable reduced Jacobian class of $f$ is simply the tranverse reduced 
Jacobian class
of $f$ acting on the transverse smooth structure, i.e. $<\phi_{u}>_{*}$.  
Similarly, the reduced  stable Jacobian class is $<\phi_{s}>_{*}$.  
It remains to see that $f$ is Anosov in the new smooth structure.  This
follows 
from Proposition 2 of Section 5, and the construction of the transverse
measure from the Gibbs measure associated to $\phi$. We postpone this
simple argument to Section 8, where explicit charts in the smooth structure
are decribed.

\begin{flushright}
$\Box$
\end{flushright}

\section{Proof of Radon-Nykodym realization}

An Anosov diffeomorphism
has a presentation as the quotient of the shift map on
a subshift of finite type \cite{Si}.  The subshift of finite type is defined
by
the transition properties of rectangles in a 
Markov partition under the action of the diffeomorphism. 
We will show that a \Ho \ cocycle over the shift map defines a ``transverse
\Ho \ measure class'' to the stable sets in the subshift of finite type.
The transverse measure class pushes down to a transverse measure class
to the stable foliation on $M$ provided the cocycle passes down
to a cocycle on $M$, i.e. when it is generated by a function on the 
subshift which is constant on fibers of the quotient to $M$.

This section is organized as follows.  In subsection 1, we recall the 
definition of a Markov partition for an Anosov diffeomorphism, and 
construct the quotient from the shift to the diffeomorphism.  
We define the notion of transverse measure class for subshifts
of finite type in subsection 2, and in subsection 3 show 
that a \Ho \ cocycle determines
a \Ho \ transverse measure class.  In subsection 4 we show that the 
transverse measure class pushes forward to the quotient space when
the cocycle does.  

\subsection{Markov partitions}

The simultaneous flow-boxes for the stable and unstable foliations of an
Anosov 
diffeomorphism $f:M \rightarrow M$ define a {\em local product structure} 
$D^{s} \times D^{u}$ on 
the manifold $M$.  A {\em rectangle} is a closed set of small diameter which
is a product in the local product structure: $R = A^{s} \times A^{u}$.  
A rectangle is {\em proper} if it is the closure of its interior.
We 
define the {\em stable} and {\em unstable boundary of R} respectively:
$\partial^{s}R = A^{s} \times \partial A^{u}$; $\partial^{u}R = \partial
A^{s} \times
A^{u}$.  

A {\em Markov partition} for $f$ is a set ${\cal C} = \lbrace R_{1}, \ldots,
R_{r}
\rbrace$ of small proper rectangles whose union is $M$, and satisfying:

\begin{description}
\item[i.] each $R_{i}$ is connected
\item[ii.] $\mbox{int}(R_{i}) \cap \mbox{int}(R_{j}) = \emptyset$ for $i
\neq j$
\item[iii.] $f(\partial^{s}{\cal C}) \subset \partial^{s}{\cal C}$ where
                              $\partial^{s}{\cal C} = \cup_{i=1}^{i=r}
\partial^{s}R$
\item[iv.] $f^{-1}(\partial^{u}{\cal C}) \subset \partial^{u}{\cal C}$ where
                        $\partial^{u}{\cal C} =
\cup_{i=1}^{i=r}\partial^{u}R$
\end{description}

Sinai proved that Anosov diffeomorphisms have Markov partitions of
arbitrarily small
diameter, where the diameter of the partition is defined to be the largest
diameter
of a rectangle in the partition \cite{Si}.

A Markov partition defines a 0-1 matrix $A$ where $A_{ij} = 1$ if 
$f(\mbox{int}R_{i}) \cap \mbox{int}R_{j} \neq \emptyset$ and is $0$ 
otherwise.
The properties of the Markov partition guarantee that if ${\bf x} \in
\Sigma_{A}$,
then the intersection $\cap_{i=-\infty}^{\infty}f^{-i}R_{x_{i}}$ consists of
a single point.  This defines a map $\pi: \Sigma_{A} \rightarrow M$ which
semi-conjugates the shift map to the mapping $f$.

\subsection{Transverse structures on a subshift of finite type.}

We recall that a subshift of finite type 
$\Sigma_{A}$ has a {\em local product structure} defined
by the stable and unstable foliations.  That is, there is a homeomorphism
defined
on a neighborhood $U$ of a point $\bf x$,
\[
u:U \rightarrow W^{s}_{\epsilon}({\bf x}) \times W^{u}_{\epsilon}({\bf x})
\]
with the property that local stable sets map to sets of the form
$W^{s}_{\epsilon}
({\bf x}) \times \lbrace {\bf z} \rbrace$, and local unstable sets map to 
sets of the form 
$\lbrace {\bf w} \rbrace \times W^{u}_{\epsilon}({\bf x})$.

A small {\em transversal} to the stable foliation
is a set which in the local product structure is represented as the graph
of a continuous function $\tau: W^{u}_{\epsilon}({\bf x}) \rightarrow 
W^{s}_{\epsilon}({\bf x})$.  

We assume that $\Sigma_{A}$ is a transitive subshift of finite type, i.e.
there is a dense orbit.  Suppose ${\bf x} \in W^{s}({\bf y})$.  Then there
is $\epsilon > 0$ and a canonical homeomorphism $h:  W^{u}_{\epsilon}({\bf
x}) 
\rightarrow W^{u}_{\epsilon}({\bf y})$, such that for every ${\bf z} \in 
W^{u}_{\epsilon}({\bf x})$, $h({\bf z}) \in W^{s}({\bf z})$.  
These homeomorphisms are called (unstable) conjugating homeomorphisms. The
{\em holonomy pseudogroup of the stable foliation} is defined to be
the pseudogroup of homeomorphisms between transversals generated by
projections
onto the unstable factor in local product charts, and unstable conjugating 
homeomorphisms between local unstable sets.

A {\em transverse measure class} to the stable foliation is an assignment of
a measure class to each small transversal, with the property that the 
holonomy transformations preserve the measure class.  A \Ho \ tranverse
measure class is one in which the representative measures on a transversal
are required to be equivalent with \Ho \ Radon-Nykodym derivative.  

A shift-invariant transverse measure class is defined as for the
diffeomorphism
case.  Note that the shift map preserves the class of transversals.  A
shift-invariant transverse \Ho \ measure class defines a \Ho \ 
Radon-Nykodym cohomology
class over the action of the shift, just as in the foliation case.

\subsection{Radon-Nykodym realization for subshifts of finite type.}

The following theorem follows easily from the standard Gibbs theory 
described in Section 5.

\begin{thm}[Radon-Nykodym realization for $\Sigma_{A}$.]
Let $\Sigma_{A}$ be a transitive subshift of finite type.  Let $<\phi>_{*}$
be
a reduced  \Ho \ cohomology class over the shift map.  Then there 
is a unique shift-invariant transverse measure class $\bmu$ to the stable
foliation 
such that the associated reduced Radon-Nykodym class is $<\phi>_{*}$.
\end{thm}

\noindent {\bf Proof.}    
We are given a \Ho \ cohomology class $<\phi>$ on $\Sigma_{A}$.  We apply 
Lemma 1 of Section 5  to obtain a \Ho \ function $\phi^{+}$ 
in the cohomology class $<\phi>$, which
we can view as a function on the one-sided shift $\Sigma_{A}^{+}$. 
We can also view $\Sigma_{A}^{+}$ as a subset of the two-sided shift, in
fact 
as a transversal to the stable foliation, as follows.  For each symbol $i
\in \lbrace
1,\ldots, r \rbrace$, pick a point ${\bf y}^{i}$ with $y^{i}_{0} = i$.  
Then $\tau = \cup _{i=1}^{r}W^{u}({\bf y^{i}},i)$ is canonically 
isomorphic to $\Sigma_{A}^{+}$,
and is a union of small transversals to the stable foliation.  Moreover,
$\tau$
meets every stable set.  We define the measure class on  
$\tau$ to be the \Ho \ measure class containing the measure
$\mu_{{\phi}^{+}}$ 
determined 
by the Bowen-Ruelle-Sinai theorem applied to $\phi^{+}$.
We define the measure class on any small transversal to be the pull-back of 
$\mu_{{\phi}^{+}}$ by 
a holonomy transformation to the transversal $\tau$ (which exists since
$\tau$ meets
every stable set).  To see that this defines  a transverse \Ho \ measure
class,
it suffices to check that the holonomy transformations between local stable
sets in
$\tau$ preserve the measure class of $\mu_{{\phi}^{+}}$, 
with \Ho \ Radon-Nykodym derivative.
But this is precisely the ``one-sided Gibbs'' property of
$\mu_{{\phi}^{+}}$.  
Finally we note that the 
Radon-Nykodym class associated to this transverse measure class is the 
reduced cohomology class of $\phi^{+} - P$, as desired.  

\begin{flushright}
$\Box$
\end{flushright}

\subsection{Pushing the transverse measure class forward.}

Let $W^{u}_{\epsilon}(x)$ be the $\epsilon$-ball about $x$ in $W^{u}(x)$. 
If $i$ is a rectangle in the Markov partition, let $W^{u}(x,i) = 
W^{u}_{\epsilon}(x) \cap R$, where $\epsilon$ is chosen so that this
is a single horizontal slice of $R$.  Let $[i]^{+} \subset \Sigma_{A}^{+}$
be
the set of sequences $y_{0}y_{1}y_{2}\ldots$ with $y_{0} = i$.  
If $x \in M$ and $x = \pi({\bf x})$ where ${\bf x} = \ldots x_{-1}x_{0}x_{1}
\ldots$ with $x_{0} =  R$, then we define a quotient map:
\[
\pi_{{\bf x},R}:[R]^{+} \rightarrow W^{u}(x,R)
\]
by $\pi_{{\bf x},R}({\bf y}) = \pi(\ldots
x_{-2}x_{-1}y_{0}y_{1}y_{2}\ldots)$.

For each rectangle $R$ and $x \in R$ with $x = \pi({\bf y})$, the measure
class on the 
transversal $W^{u}(x,R)$ is defined 
to be the image by $\pi_{{\bf y},R}$ of the measure class $\bmu_{\phi}$ on 
$[R]^{+}$.  

We need to check that if parts of $W^{u}(x,R)$ and
$W^{u}(x^{\prime},R^{\prime})$
correspond under local projection along the stable foliation, then the
measure
classes defined by $\pi_{{\bf x},R}$ and $\pi_{{\bf x}^{\prime},R^{\prime}}$
also correspond.  If $R = R^{\prime}$, this follows
from the fact that, if $p_{x,x^{\prime}}$ is the projection along local
stable
leaves 
from $W^{u}_{\epsilon}(x)$ to $W^{u}_{\epsilon}(x^{\prime})$, and if
$\pi({\bf x}) = x$
and $\pi({\bf x}^{\prime}) = x^{\prime}$, then 
\[
\pi_{{\bf x},R} =  p_{x,x^{\prime}} \circ \pi_{{\bf x}^{\prime},R}
\]
If $R \neq R^{\prime}$ we consider representative measures on the
$W^{u}(x,R)$
and $W^{u}(x^{\prime}, R^{\prime})$.  Let $\mu_{x,R}$ denote the image of 
$\mu_{{\phi}^{+}}$ restricted to $[R]^{+}$ by $\pi_{{\bf x},R}$.

Let $U \subset W^{u}(x,R)$ and 
$V \subset W^{u}(x^{\prime},R^{\prime})$ be such that $p_{x,x^{\prime}}:
U \rightarrow V$ is a homeomorphism.  
There are two main points:

\begin{description}
\item[1.   ] Let  $y \in U$, and  $y^{\prime} = p_{x,x^{\prime}}(y)$.  
We can make sense of the expression $\Phi^{+}(y \rightarrow y^{\prime})$ 
{\em on $U$} as follows.  Recall
\[
\Phi(y \rightarrow y^{\prime}) = \Sigma_{k=0}^{\infty}(\phi \circ
f^{k}(y^{\prime})
- \phi \circ f^{k}(y)).
\]
is defined and \Ho \ on $U$.
Note that the transfer function $u:\Sigma_{A} \rightarrow {\bf R}$ which
makes
$\phi$ cohomologous to $\phi^{+}$ pushes forward to a well-defined and \Ho \
function $u_{R}$ on each of the individual quotients 
$W^{u}(x,R)$   Therefore we {\em define} a \Ho \ function:
\[
\Phi^{+}(y \rightarrow y^{\prime}) = \Phi(y \rightarrow y^{\prime}) + 
u_{R^{\prime}}(y^{\prime})
- u_{R}(y).
\]

\item[2.   ] Let $\mu_{U}$ be $\mu_{{\bf x},R}$ restricted to $U$.  Define
$\mu_{V}$
similarly.   Then
\[
\frac{d((p_{x,x^{\prime}})_{*}(\mu_{V}))}{d\mu_{U}} = {\rm exp}(\Phi^{+}(y
\rightarrow y^{\prime})).
\]
\end{description}

\noindent We need to prove the second statement.  We define a measure $\nu$
supported
on $V$ by $\nu = p_{x,x^{\prime}}^{*}({\rm exp}(\Phi^{+})\mu_{U})$.  
We want to pull back to $\Sigma_{A}^{+}$
both $\nu$ and $\mu_{V}$, where we will show them to be equal.
For this to make sense, and imply the second statement, we need the
following lemma.
Let $\tilde{V} = \pi_{x,R^{\prime}}^{-1}(V)$.

\begin{lem}  
$\pi_{{\bf x},R^{\prime}}:\tilde{V} \rightarrow V$ is one-to-one over a set
of full
measure with respect to both $(p_{x,x^{\prime}})^{*}\mu_{U}$ and $\mu_{V}$.
\end{lem}

\noindent {\bf Proof.}  Let $Y \subset V$ be the set of points with more
than one
preimage by $\pi_{{\bf x},R^{\prime}}$.  We need to show that 
$\mu_{U}(p_{x,x^{\prime}}^{-1}(Y)) = 0$ and $\mu_{V}(Y) = 0$

Let $\mu$ be the invariant Gibbs measure on $M$ associated to the cohomology
class
$<\phi>$. We need the following two facts.  The technical proofs are
included in 
the appendix.

\begin{description}

\item[1.   ] Let $R = A^{s} \times A^{u}$ in the local product structure.
Let 
$x \in R$, and $Y \subset W^{u}(x,R)$.  Then $\mu(A^{s} \times Y) = 0$
implies
$\mu_{x,R}(Y) = 0$. 

\item[2.   ] $\mu(\partial {\cal C}) = 0$.

\end{description}

\noindent Now we prove the lemma.  If $y \in Y$, then $d(f^{n}(y),\partial
{\cal C})
\rightarrow 0$ as $n \rightarrow \infty$.  The same is therefore true of
points
in $Y^{\prime} := p_{x,x^{\prime}}^{-1}(Y)$. 
 Let $R = A^{s}(R) \times A^{u}(R)$ and $R^{\prime}
= A^{s}(R^{\prime}) \times A^{u}(R^{\prime})$.  Let $Z = A^{s}(R) \times Y$
and $Z^{\prime} = A^{s}(R^{\prime}) \times Y^{\prime}$.  Then
$d(f^{n}(Z),\partial {\cal C}) \rightarrow 0$ and
$d(f^{n}(Z^{\prime},\partial
{\cal C}) \rightarrow 0$ as $n \rightarrow \infty$.  Since $\mu(\partial
{\cal C}) = 0$
, and $\mu$ is invariant, 
$\mu(Z) = 0$ and $\mu(Z^{\prime}) = 0$, and we conclude that 
$\mu_{V}(Y) = 0$ and $\mu_{U}(Y^{\prime}) = 0$.

\begin{flushright}
$\Box$
\end{flushright}

We return to the proof of statement 2 relating the measures on $U$ and $V$.
We want to show that the measure $\tilde{\nu} = (\pi_{x,R^{\prime}})_{*}\nu$
coincides with $(\pi_{x,R^{\prime}})_{*}\mu_{V}$.  The latter is simply
$\mu_{\phi^{+}}$ restricted to $\tilde{V}$.  By the Local Uniqueness
corollary 
to the Bowen-Ruelle-Sinai theorem, it suffices to show that $\tilde{\nu}$
has Radon-Nykodym derivative ${\rm exp}(-\sum_{k=0}^{n-1}\phi^{+}\circ 
\sigma^{k} + nP)$ under ${\bf x} \rightarrow \sigma^{n}({\bf x})$ whenever
${\bf x}$ and $\sigma^{n}({\bf x})$ are both in $\tilde{V}$.

Recall $\nu = p_{x,x^{\prime}}^{*}({\rm exp}(\Phi^{+})\mu_{U})$.  
So 
\[
\frac{d({\tilde{\nu}(\sigma({\bf y})}))}{d\nu({\bf y})} =
\pi^{*}(p_{x,x^{\prime}})_{*}
\frac{d({\rm exp} (\Phi^{+})\mu_{U})(f(y))}{d({\rm
exp}(\Phi^{+})\mu_{U})(y)}
\]
\noindent The main point in the calculation is the following.  $\phi^{+}$ is
well-defined
on the individual quotients $W^{u}(x,R)$, namely $\phi^{+}_{R} = \phi +
u_{R} \circ
f - u_{R}$.  Then 
\[
\Phi^{+}(f(y)\rightarrow f(y^{\prime}))
- \Phi^{+}(y \rightarrow y^{\prime}) = -\phi^{+}_{R^{\prime}}(y^{\prime}) + 
\phi^{+}_{R}(y)
\]
where
$y^{\prime} = p_{x,x^{\prime}}(y)$ and $f(y^{\prime}) = f(p_{x,x^{\prime}})
= p_{x,x^{\prime}}(f(y))$.  

Thus 
\begin{eqnarray}
        \frac{d({\rm exp} (\Phi^{+})\mu_{U})(f(y))}{d({\rm
exp}(\Phi^{+})\mu_{U})(y)}
                & = & {\rm exp}(-\phi^{+}_{R^{\prime}}(y^{\prime}) + 
                        \phi^{+}_{R}(y)) \cdot {\rm exp} (-\phi^{+})\\
                & = & {\rm exp}(-\phi^{+}_{R^{\prime}}(p_{x,x^\prime}(y)) +
P)
\end{eqnarray}
which is the desired result.

\section{Gibbs charts}

A Markov partition for a hyperbolic automorphism $L$ of $T^{2}$ can be
constructed as
follows \cite{AW}.  Let $E^{s}$ and $E^{u}$ be the stable and 
unstable eigenspaces respectively.
Project into the torus a segment in $E^{s}$ through the origin, and a
segment in $E^{u}$
through the origin.  Extend these segments until they cut the torus into
parallelograms.
The segment in the stable direction should map into itself under $L$, and
the segment in
the unstable direction should map into itself under $L^{-1}$.   This
decomposition 
of the torus is a Markov partition.  Let $A$ be the transition matrix of the
partition,
and let $\pi:\Sigma_{A} \rightarrow T^{2}$ be the quotient map from the
subshift of 
finite type defined by $A$.  The unstable segment $\tau_{u}$ 
(which is also the unstable boundary
of the partition) is the image by $\pi$ of a  copy of the {\em one-sided
shift} $\Sigma_{A}^{+}$ specified by fixing the backward part of a sequence
to be 
the backward part of some fixed pre-image of the origin.  Similarly, the 
stable segment $\tau_{s}$ is the image of a copy of the backward one-sided
shift defined
by $A$, or equivalently of the one-sided shift defined by the $A^{t}$, the
transpose of $A$.

The smooth structure determined by a pair of  reduced cohomology classes
 $<\phi_{u}>_{*}$ and $<\phi_{s}>_{*}$  has the following explicit
description.
We can assume that the functions $\phi_{u}$ and $\phi_{s}$ have pressure 0.
If not,
we can just subtract the pressure, which will not change the reduced
cohomology class.
Let $\tilde{\phi_{u}}$ be the pull-back of $\phi_{u}$ to $\Sigma_{A}$, by
$\pi$.
$\tilde{\phi_{s}}$ is defined similarly.
Now change $\tilde{\phi_{u}}$ by a coboundary to get a function
$\phi_{u}^{+}$ that 
depends only on the forward part of a sequence.  Similarly, one gets a
function
$\phi_{s}^{-}$ which is cohomologous to $\tilde{\phi_{s}}$ and depends only
on the 
backward part of a sequence.  We can regard these as functions on
$\Sigma_{A}^{+}$
and $\Sigma_{A^{t}}^{+}$, respectively.
Now let $\mu_{+}$ be the unique probability measure on $\Sigma_{A}^{+}$
satisfying
\[
\frac{d\mu_{+}(\sigma({\bf x}))}{d\mu_{+}({\bf x})} = {\rm
exp}(\phi_{u}^{+}({\bf x}))
\]
Let $\mu_{-}$ be the unique probability measure on $\Sigma_{A^{t}}^{+}$
with 
\[
\frac{d\mu_{-}(\sigma({\bf y}))}{d\mu_{-}({\bf y})} = {\rm
exp}(\phi_{s}^{-}({\bf y}))
\]
These measures push-forward to measures $\nu^{+}$ and $\nu^{-}$ on
$\tau_{u}$ and 
$\tau_{s}$ respectively. 

Smooth coordinate charts (in the new structure) are obtained by integrating
the measures
$\nu^{+}$ and $\nu^{-}$ along two side of a "rectangle" in the partition,
and 
taking the product structure.   The main point in the proof of the theorem
is that 
if the segment $\tau_{u}$ is presented {\em differently} as the image of
$\Sigma_{A}^{+}$,
i.e. by choosing a different pre-image of the origin, with different
backward part,
then the smooth coordinate along $\tau_{u}$ determined by pushing forward
the measure
$\nu^{+}$ by this different presentation is smoothly equivalent to the
original one.

It is now clear from Proposition 2 of Section 5 that $f$ is Anosov in the
new 
smooth structure.  We simply use the expanding measure in the \Ho \ measure
class
associated to $<\phi_{u}>_{*}$ to define the coordinate along $\tau_{u}$.
Similarly,
use the expanding measure associated to $<\phi_{s}>_{*}$ to define the
coordinate 
along $\tau_{s}$.  

\noindent {\bf Remark.}  The two maps $\pi_{1}:\Sigma_{A}^{+} \rightarrow
\tau_{u}$ and
$\pi_{2}:\Sigma_{A}^{+} \rightarrow \tau_{u}$, determined by viewing
$\tau_{u}$ 
from the ``clockwise'' side or the ``counterclockwise'' side, determine a
comparison
map $\alpha:\Sigma_{A}^{+} \rightarrow \Sigma_{A}^{+}$ on a set of full
measure,
defined by $\pi_{2}(\alpha({\bf x})) = \pi_{1}({\bf x})$. 
The measure $\mu^{+}$ has \Ho \ Radon-Nykodym derivative under $\alpha$.

\newpage

\appendix

\begin{flushleft} 
{\bf \LARGE Appendix}
\end{flushleft}

\bigskip

\section{Quasisymmetric equivalence of Gibbs structures}

We show that the conjugacy  between the linear map
in a topological conjugacy class, and the map constructed from  a \Ho \ 
cocycle, is quasisymmetric along the leaves of the stable and unstable
foliations.   The proof is adapted from \cite{Ji}.

A homeomorphism $h: I \rightarrow I$ is {\em quasisymmetric} if
there exists a $K > 0$ such that for every pair of adjacent intervals
$I$ and $J$ of equal length,
\[
1/K \leq \frac{|h(I)|}{|h(J)|} \leq K.
\]
The number $K$ is the {\em quasisymmetry constant} of 
the map, and $h$ is
called $K-quasisymmetric$.
The composition of a quasisymmetric map with a $C^{1}$ map is again
quasisymmetric.

We recall Proposition 2 from Section 5.  Let $\Sigma_{A}^{+}$ be a mixing
one-sided subshift of finite type.  Let $\pi: \Sigma_{A}^{+} \rightarrow
I$ be a \Ho \ map onto an interval $I \subset {\bf R}$, with the property
that $\pi(\sigma({\bf x})) = \pi(\sigma({\bf y}))$ whenever $\pi({\bf x})
 = \pi({\bf y})$.  In addition we assume that the image by $\pi$ of 
a cylinder set $C_{n}$ is an interval.  
Let $\phi^{\prime}$ be a \Ho \ function on $I$, and let 
$\phi$ be the pull-back  of $\phi^{\prime}$ to $\Sigma_{A}^{+}$ by $\pi$.
Let $\mu_{0}$ be the measure of maximal entropy on $\Sigma_{A}^{+}$,
and let $\mu_{\phi}$ be the measure associated to $\phi$ as constructed in
the 
Bowen-Ruelle-Sinai theorem, that is, the unique probability measure with
Radon-Nykodym derivative $\phi - P(\phi)$ where $P(\phi)$ is the pressure
of $\phi$.  Assume that $\pi$ is injective on a set of full measure, with
respect to both $\mu_{0}$ and $\mu_{\phi}$.   There are two metrics $d_{0}$
and
$d_{\phi}$ on $I$, defined by the push-forward by $\pi$ of $\mu_{0}$ and
$\mu_{\phi}$,
respectively.

\begin{prop}[Quasisymmetric equivalence of Gibbs structures]
The identity map $\iota:(I,d_{0}) \rightarrow (I,d)$ is quasisymmetric.
\end{prop}

We note that the hypotheses of the proposition are 
satisfied by the map $\pi:\Sigma_{A}^{+} \rightarrow \tau_{u}$ defined
in section 8.  

The proof given in \cite{Ji} applies in 
a more general context.  The basic idea is the following.  The image
of the partitions into cylinder sets of the subshift of finite type
defines a nested sequence of partitions of the interval $I$.  The identity
map of course preserves these partitions.  The quasisymmetry estimate
follows from bounds on the geometry of this sequence of partitions, in both
the $d_{0}$ and the $d_{\phi}$ metrics.

We can assume that $\phi$ is in fact the Radon-Nykodym derivative of the
expanding (equilibrium) measure on $\Sigma_{A}^{+}$ 
associated to the reduced cohomology class 
$<\phi>_{*}$. This is because the expanding measure, and the measure
$\mu_{\phi}$ are equivalent with a \Ho \ Radon-Nykodym derivative,  
and therefore the corresponding metrics on $I$ are $C^{1 + \alpha}$
equivalent.

Let ${\cal C}_{n}$ denote the partition of $\Sigma_{A}^{+}$
by cylinder sets of size $n$.
Let ${\cal D}_{n}$ be the corresponding sequence of partitions of $I$.

In the following lemma, we consider $I$ with the $d_{0}$ metric, and show
how to approximate the intervals defined by an equally spaced triple of 
points by 
elements of the partitions ${\cal D}_{n}$.  

\begin{lem}
There exists a positive integer $N = N(\Sigma_{A}^{+})$ with the following
property.
Let $x,y \in I$, and let $z$ be the midpoint of the interval $[x,y]$.  
We will write $R = [x,z]$ and $S = [z,y]$.  Let 
$n$ be the smallest integer such that there exists $D \in {\cal D}_{n}$ with
$R \cup S \subset D$.  Then there are $D_{R}, D_{S} \in {\cal D}_{n + N}$
with $D_{R} \subset R$ and $D_{S} \subset S$.
\end{lem}

\noindent {\bf Outline of the proof.}  Let $h$ be the topological entropy of
$\Sigma_{A}^{+}$,
and let $M$ be the mixing time, that is, $A^{M}$ has all positive entries.  
Then $N \geq 4M + 1 + \frac{log 2}{log h}$ has
this property.  This follows easily from the following three geometric
properties
of the partitions ${\cal D}_{n}$ in the $d_{0}$ metric.  Let $\lambda = \exp
h$.
We denote the length of an interval $D$ in the $d_{0}$ metric by $|D|_{0}$.
\begin{description}

\item[1.  Exponentially decreasing geometry.]   

Let $D_{n + m} \in {\cal D}_{n + m}$, $D_{n} \ in {\cal D}_{np}$, with 
$D_{n + m} \subset D_{n}$.  Then
\[
\frac{|D_{n + m}|_{0}}{|D_{n}|_{0}} \leq \lambda^{m - M}.
\]

\item[2.  Bounded ratio geometry.]

Let $D_{n} \in {\cal D}_{n}$, $D_{n + 1} \in {\cal D}_{n + 1}$, with\break
$D_{n + 1} \subset D_{n}$.  Then
\[
\frac{|D_{n + 1}|_{0}}{|D_{n}|_{0}} \geq \lambda^{-(M + 1)}.
\]

\item[3.  Bounded nearby geometry.]

Let $D, E \in {\cal D}_{n}$ be adjacent intervals.  Then
\[
\lambda^{-M} \leq \frac{|D|_{0}}{|E|_{0}} \leq \lambda^{M}.
\]
\end{description}

\begin{flushright}
$\Box$
\end{flushright}

\noindent {\bf Outline of proof of the proposition.}
Analagous geometric properties hold for the partitions ${\cal D}_{n}$ in the
$d_{\phi}$ metric.  The quasisymmetry estimate follows from these.
We define some preliminary quantities.

\bigskip

\noindent Let ${\rm S}_{n}\phi ( {\bf x}) = \sum_{k = 0}^{n - 1} 
(\phi(\sigma^{k}({\bf x}))$.
Let 
\[{\rm var}_{n} \phi = {\rm sup} \lbrace | {\rm S}_{n}\phi({\bf x}) - 
{\rm S}_{n}\phi({\bf y}) | 
\ {\rm where} \  {\bf x},{\bf y} \in C_{n} \  {\rm for \ some} \  
C_{n} \in {\cal C}_{n} \rbrace.
\]
Since $\phi$ is \Ho, there exist $c > 0$ and $\beta < 1$ such that 
${\rm var}_{n} \phi < c \beta^{n}$.  Therefore 
\[
{\rm var} \phi =: \sum_{k= 0}^{\infty} {\rm var}_{k} \phi < \infty.
\]
Let $\| \phi \| = {\rm sup} \lbrace | \phi({\bf x}) | : {\bf x} \in
\Sigma_{A}^{+} \rbrace$.
Define  $L = \exp (2 {\rm var} \phi) \cdot \exp ( M \| \phi \| )$.

\bigskip

\noindent In what follows, all lengths are with respect to the $d_{\phi}$
metric,
where $\phi$ is the Radon-Nykodym derivative of the expanding measure 
on $\Sigma_{A}^{+}$
associated to the reduced cohomology class $<\phi>_{*}$.  The $d_{\phi}$
length
of an interval $D$ is denoted $|D|_{\phi}$.
\begin{description}

\item[1.  Bounded ratio geometry.]
Let $D_{n + 1} \in {\cal D}_{n + 1}$, $D^{n} \in {\cal D}_{n}$, with\break
$D_{n + 1} \subset D_{n}$.  Then
\[
\frac{|D_{n + 1}|_{\phi}}{|D_{n}|_{\phi}} \geq L^{-1} \cdot \exp (- \| \phi
\|) .
\]

\item[2.  Bounded nearby geometry.]
Let $D, E \in {\cal D}_{n}$ be adjacent intervals.  Then
\[
\frac{1}{L} \leq \frac{|D|_{\phi}}{|E|_{\phi}} \leq L.
\]

\end{description}

\noindent These properties follow from Bowen's estimate for the 
$\mu_{\phi} - $measure
of a cylinder set $C_{n}$, when $\phi$ is the Radon-Nykodym derivative of
the invariant measure.\cite{B}  For any ${\bf x} \in C_{n}$,
\[
c_{1} \leq \frac{\mu_{\phi}(C_{n})}{exp (S_{n}\phi({\bf x}))} \leq c_{2}
\]
where $c_{1} = \exp (-M \| \phi \| ) \cdot \exp (-{\rm var} \phi)$, and
$c_{2} = \exp ( {\rm var}\phi)$. 

Using these properties, and the Lemma, one obtains an estimate for the
quasisymmetry constant $K$:
\[
K \leq L^{(2N + 2)} \cdot \exp ((2N + 1) \|  \phi \|).
\]

\noindent {\bf Remark.}  This estimate is not sharp, as can be seen by 
considering $\phi$ to be the constant function with pressure zero.  
The estimate can be improved by a more careful comparison of the
partition elements lying in the pair of adjacent intervals.

\section{Gibbs measures and Markov partition boundaries}

We prove here the technical facts needed in Section 7.4  Let $f:M
\rightarrow M$ be
an Anosov diffeomorphism, $\cal C$ be a Markov partition for $f$, and $\mu$
be the 
Gibbs measure associated to the cohomology class $< \phi >$. 

\begin{prop}
Let $R \in {\cal C}$ be a rectangle, with $R = A_s \times A_u$ in the local 
product structure.  Let $x \in R$, and $Y \subset W^{u}(x,R)$.  Let
$\mu_{x,R}$ 
be the push-forward of the one-sided measure $\mu_{\phi^{+}}$,
as defined in Section 7.4.
Then \hbox{$\mu(A^{s} \times Y) = 0$} implies $\mu_{x,R}(Y) = 0$.
\end{prop}

\noindent {\bf Proof.} 
We recall how $\mu$ on $T^{2}$ can be constructed from 
$\mu_{\phi^{+}}$ on $\Sigma_{A}^{+}$. See \cite{B} for details.
If $\pi: \Sigma_{A} \rightarrow T^{2}$
is the quotient map determined by the Markov partition, then $\mu$ is the 
push-forward by $\pi$ of the measure $\tilde{\mu}$ on $\Sigma_{A}$, defined
as follows.  Let $\nu$ be the {\em shift-invariant} probability measure on 
$\Sigma_{A}^{+}$ equivalent to $\mu_{\phi^{+}}$. (The Radon-Nykodym
derivative
of $\nu$ with respect to $\mu_{\phi^{+}}$ is given explicitly up to a
constant
factor as the unique positive 
eigenvector of the Perron-Fr\"{o}benius operator associated to $\phi^{+}$.)
The measure $\tilde{\mu}$ on $\Sigma_{A}$ is obtained from the measure 
$\nu$ on $\Sigma_{A}^{+}$ as follows.  If $g$ is a continuous function on
$\Sigma_{A}$,
define $g^{*}$ on $\Sigma_{A}^{+}$ by 
\[
g^{*}({\bf x}) = {\rm min} \lbrace g({\bf y})\ {\rm where} \ y_{i} = x_{i} \
{\rm for \  all} \  i \geq 0 \rbrace
\]
Then $\lim_{n \rightarrow \infty} \nu ( (g \circ \sigma^{n})^{*})$ exists,
and we
define $\tilde{\mu}(g)$ to be this limit. This linear functional defines the
measure $\tilde{\mu}$.  The proposition now follows easily.

\begin{flushright}
$\Box$
\end{flushright}

\begin{prop} 
$\mu(\partial{\cal C}) = 0$.
\end{prop}

\noindent  {\bf Proof.}  The following proof is adapted from Bowen's in the
case $\mu$
is the measure associated to the constant cocycle. \cite{B1}  The proof
relies on the
{\em variational principle}, and the fact that the topological
pressure
always decreases when the dynamics is restricted to an invariant subset.

\noindent {\bf Variational Principle.}  
If a  cohomology class over $f$, say $< \phi >$ has been fixed, we define
the 
{\em measure theoretic pressure} of an invariant probability measure $\nu$
to 
be $h_{\nu} + \int \phi d\nu$, where $h_{\nu}$  is the measure theoretic
entropy
of $f$.  We denote this $P_{\nu}(\phi,f)$.  The topological pressure of
$\phi$  
is denoted $P(\phi,f)$.  Note that these depend only on the cohomology class
of 
$\phi$.  The variational principle states:
\[
P_{\nu}(\phi,f) \leq P(\phi,f)  \ {\rm and } \ \sup_{\nu} P_{\nu}(\phi,f) =
P(\phi,f)
\]
where the supremum is over all invariant probability measures.  The
variational 
principle is true for any homeomorphism of a compact metric space.  If the
mapping
is Anosov, and the cocycle is \Ho \, then 
supremum is achieved precisely at the Gibbs measure $\mu$ associated to
$<\phi>$.
Recall that if $W$ is a compact
f-invariant subset, then $P(\phi_{\mid W}, f_{\mid W}) < P(\phi,f)$.

Let $\mu$ be the Gibbs measure associated to the cohomology class $<\phi>$.
We will show that $\mu(\partial^{s}{\cal C}) = 0$.  A similar argument shows
that 
$\mu(\partial^{u}{\cal C}) = 0$.  Let $W = \cap_{n \geq 0}
f^{n}(\partial^{s}{\cal C})$.
Suppose that $\mu(\partial^{s}{\cal C}) = a > 0.$  Then $\mu(W) = a$, and 
$\mu(\cup f^{n}(\partial {\cal C}) \setminus W) = 0$.  Define $\nu_{1}$ on
$W$ by 
$\nu_{1} = \frac{1}{a} \mu$, and $\nu_{2}$ on $M$ by $\nu_{2}(E) =
\frac{1}{1-a}
\mu(E \setminus W)$.  Then $\nu_{1}$ and $\nu_{2}$ are f-invariant, have 
disjoint support, and $\mu = a \nu_{1} + (1 - a) \nu_{2}$.  Therefore
\[
P_{\mu}(\phi,f) = a P_{\nu_{1}}(\phi,f) + (1 - a) P_{\nu_{2}}(\phi,f).
\]
The variational principle implies that $P_{\nu_{2}}(\phi,f) \leq P(\phi,f)$
and 
\begin{eqnarray}
P_{\nu_{1}}(\phi,f) &  = & P_{\nu_{1}}(\phi_{\mid W}, f_{\mid W}) \\
		    &  \leq & P(\phi_{\mid W},f_{\mid W}) \\
                    &  < & P(\phi,f).  
\end{eqnarray}
But then $P_{\mu}(\phi, f) < P(\phi,f)$, a contradiction
since $\mu$ achieves the supremum.

\begin{flushright}
$\Box$
\end{flushright}
\bigskip\bigskip

\bibliographystyle{alpha}
\bibliography{ref}

\end{document}